\documentclass[12pt]{amsart}

\usepackage{amsmath,amsthm,amsfonts,latexsym,amssymb,enumerate,color}
\usepackage{graphicx}
\usepackage{enumitem}  
   
\def\FF{\mathcal F}
\newcommand\CC{{\mathbb C}}

\newcommand\GF{{\mathcal F}}
\newcommand\GD{{\mathcal D}}
\newcommand\DD{{\mathbb D}}

\newcommand\NN{{\mathbb N}}

\newcommand\RR{{\mathbb R}}

\newcommand\TT{{\mathbb T}}

\newcommand\ZZ{{\mathbb Z}}

\newcommand\Id{{\rm Id}}

\def\beq{\begin{equation}}
\def\eeq{\end{equation}}
\newtheorem{thm}{Theorem}[section]
\newtheorem{prop}[thm]{Proposition}

\newtheorem{lem}[thm]{Lemma}
\newtheorem{cor}[thm]{Corollary}
\newtheorem{rem}[thm]{Remark}

\newcommand\beginpf{\noindent {\bf Proof:} \quad}

\newcommand\re{\mathop{\rm Re}\nolimits}

\newcommand\esssup{\mathop{\rm ess\ sup}\nolimits}

\def\beginpf{\begin{proof}}
\def\endpf{\end{proof}}

\renewcommand\phi{\varphi}

\DeclareMathOperator{\dom}{dom}

\def\Re{\re}

\newcommand\Hol{\mathop{\rm Hol}}
\newcommand\LL{\mathcal L}
\newcommand\Vco{W_{v,\psi}}
\newcommand\W{W_{w,\phi}}

\begin{document}
\title[Weighted composition operators]{Weighted composition operators on the Fock space: iteration and semigroups}

\author{I. Chalendar}
\address{Isabelle Chalendar,  Universit\'e Gustave Eiffel, LAMA, (UMR 8050), UPEM, UPEC, CNRS, F-77454, Marne-la-Vallée (France)}
\email{isabelle.chalendar@univ-eiffel.fr}

\author{J.R. Partington}
\address{Jonathan R. Partington, School of Mathematics, University of Leeds, Leeds LS2 9JT, UK}
 \email{J.R.Partington@leeds.ac.uk}
 
 
\baselineskip18pt

\bibliographystyle{plain}

\begin{abstract}
This paper considers discrete and continuous semigroups of (weighted) composition
operators on the Fock space. For discrete semigroups consisting of powers of a single operator, the
asymptotic behaviour of the semigroups is analysed. For continuous semigroups and groups, 
a full classification of possible semigroups is given, and the generator is calculated.
\end{abstract}

 \maketitle

\section{Introduction and notation}

For $1 \le \nu < \infty$ the Fock space $\GF^\nu$ consists of all entire functions $f: \CC \to \CC$ such that the norm
\[
\|f\|_{\nu}:= \left( \frac{ \nu}{2\pi}\int_\CC |f(z)|^\nu e^{- \nu |z|^2/2} \, dm(z) \right)^{1/\nu} 
\]
is finite. 

Note that for $k \ge 0$ we have
\[ 
\int_0^\infty r^k e^{-r^2} dr = \frac{1}{2} \Gamma\left(\frac{k+1}{2}\right),
\]
from which we see that
the functions $e_n$, $n=0,1,2,\ldots$  such that
\beq\label{eq:onb}
e_n(z) =  z^n / \sqrt{n!}
\eeq
form an orthonormal basis for the Hilbert space $\GF^2$.

For entire functions $w,\phi,f$ we define the composition operator 
$C_\phi$ by $C_\phi f=f \circ \phi$ and the weighted composition operator $\W$ by $\W f = w C_\phi f$.

Carswell, MacCluer and Schuster \cite{CMS03} characterised the bounded composition operators 
$C_\phi$
on $\GF^2$ as follows:
\begin{prop}\label{prop:11}
The operator $C_\phi$ is bounded on $\GF^2$ if and only if one of the following conditions holds:
\begin{itemize}
\item $\phi(z)= az+b$ with $|a|<1$, $b \in \CC$;
\item $\phi(z)=az$ with $|a|=1$.
\end{itemize}
In the case $|a|<1$ the operator $C_\phi$ is compact.
\end{prop}

For weighted composition operators $\W$ we have the result of Le \cite{Le}, extended by
Hai and Khoi \cite{HaiKhoi} as follows:

\begin{prop}\label{prop:wcobdd}
The weighted composition operator $\W$ with $w$ not identically zero is bounded on $\GF^\nu$ if and only if both the following
conditions hold:
\begin{itemize} 
\item $w \in \GF^\nu$;
\item $M(w,\phi):=\sup_{z \in \CC} |w(z)|^2 e^{(|\phi(z)|^2-|z|^2)} < \infty$.
\end{itemize}
Moreover, in this case $\phi(z)=az+b$ with $|a| \le 1$. If $|a|=1$ then we also
have 
\beq\label{eq:wamod1}
w(z)=w(0)e^{- \overline b az}
\eeq
 for $z \in \CC$.
\end{prop}

See also \cite{CG21, Men21, Men22}  for a recent discussion of these results.\\

Zhao \cite{zhao} gave a criterion for the invertibility of a weighted composition operator
on the Fock space; using the more recent results above we can simplify this as follows.   

\begin{prop} \label{prop:invert}
A weighted composition $\W$ is invertible on $\GF^\nu$ if and only if $\phi(z)=az+b$ with
$|a|=1$, and $w(z)=w(0)e^{- \overline b az}$ with $w(0) \ne 0$.
\end{prop}
\beginpf
Suppose that $\W$ is invertible. Then $w$ must never vanish as if it did then any zeros of $w$ would also be zeros of
$\W f$ for every $f \in \GF^\nu$ and $\W$ would not be surjective. For $f \in \GF^\nu$ set
$g(z)=w(z)f(az+b)$ with $0<|a| \le 1$. Then $f(\zeta)w(a^{-1}\zeta - a^{-1}b)=g(a^{-1}\zeta - a^{-1}b)$
for $\zeta \in \CC$. Hence the inverse of $\W$ is also a weighted composition operator
$\Vco$ with $\psi(\zeta)=\phi^{-1}(\zeta)=a^{-1}\zeta - a^{-1}b$, and thus $|a|=1$ if $\Vco$  is bounded.
We also have a formula for $w$, namely \eqref{eq:wamod1}.
Similarly,
\[
v(z)=v(0)e^{a\overline b a^{-1}z}=v(0)e^{\overline b z}
\]
and since
\[
1/w(a^{-1}z - a^{-1}b) = w(0)^{-1} e^{\overline b a (a^{-1}z - a^{-1}b) }
\]
we find that $v(z)= 1/w(a^{-1}z - a^{-1}b)$ if $v(0)=w(0)^{-1}e^{-|b|^2}$.
\endpf

This result will be of 
further use to us when we discuss groups of weighted composition operators.\\

The paper is organised as follows. In Section~\ref{sec:2} we study the asymptotics of iteration of bounded  composition operators on Fock spaces. Theorem~\ref{th:asymp} is a complete answer to this problem. 
Then we consider $C_0$ semigroups of composition operators $(C_{\varphi_t})_{t\geq 0}$, providing explicit expressions of the semiflow $(\varphi_t)_{t\geq 0}$ in Theorem~\ref{thm:compsg}. We also prove that 
every $C_0$ semigroup $(T_t)_{t\geq 0}$ on ${\mathcal F}^2$ with generator  $Af=Gf'$ is a semigroup of composition operators,
which implies that 
$G(z)=az+b$ with $\Re (a)\leq 0$ and $b\in\CC$.  Moreover the convergence in norm as $t\to\infty$ of such $(T_t)_{t\geq 0}$ is equivalent to  the condition $\Re (a)<0$. We end this section with 
Theorem~\ref{thm:analytic}, which is a criterion, on  the space ${\mathcal F}^2$ and for a generator $A$, for generating an analytic semigroup of composition operators.  

Section~\ref{sec:3} is concerned with weighted composition operators $W_{w,\varphi}$. We first characterize, in Theorem~\ref{thm:3.1}, the weights for which $ W_{w,\varphi}$ is power-bounded 
on ${\mathcal F}^\nu$ in the particular case where $\varphi(z)=az+b$ with $|a|=1$. The case $|a|<1$ is much more delicate. We restrict our study to the case when the weight $w$ does not vanish, which is the relevant situation for semigroups of weighted composition operators, as shown in Theorem~\ref{thm:c0sg}.    We consider two cases and split our results into Theorem~\ref{thm:3.3} and Theorem~\ref{thm:3.4}.  We also derive a result on the converegnce on the asymptotic behaviour of the iterates of $W_{w,\varphi}$ in Corollary~\ref{cor:cv}, thanks to the description of its eigenvalues and eigenspaces described in Theorem~\ref{th:eigen}. 
The last part of the paper is a complete description of $C_0$ semigroups and $C_0$ groups of weighted composition operators on ${\mathcal F}^\nu$.


\section{Composition operators}\label{sec:2}

We begin with (unweighted) composition operators, where some fairly complete results will
be derived.

\subsection{Iteration}

Note first that for a composition operator $C_\phi$ with $\phi(z)=az+b$, and $a \ne 1$ we have
\[
C_\phi^n f(z)=f\left(a^n z + \frac{1-a^n}{1-a}b\right) \qquad \hbox{for} \quad n=1,2,3,\ldots,
\]
as is easily seen by induction. Also, for $\phi(z)=z+b$ we clearly have
\[
C_\phi^n f(z)=f(z+nb) \qquad \hbox{for} \quad n=1,2,3,\ldots.
\]

There is a general convergence criterion for  powers of an operator $T \in \LL(X)$,
where $X$ is a Banach space, which may be found in \cite[A-III, Sec. 3.7]{crowd}; it was
applied to composition operators in \cite{ACKS}, and to weighted composition
operators in \cite{CP21}.

An operator $T \in \LL(X)$ is power-bounded if $\sup_{n\geq 1}\|T^n\|<\infty$. 
 
\begin{prop}\label{prop:2.1}
Let $T \in \LL(X)$  be power-bounded.  Then the following assertions are equivalent:\\
(i) $P:=\lim T^n$ exists in $\LL(X)$ and is of finite rank.\\
(ii) (a) $r_e(T)<1$;\\
(b) $\sigma_p(T) \cap \TT \subseteq \{1\}$;\\
(c) if $1$ is in $\sigma(T)$ then it is a pole of the resolvent of order 1.\\
In that case, $P$ is the residue at $1$.
\end{prop}

The following formula for $\|C_\phi\|$ is given in \cite{CMS03} (see also \cite{dai}).
If $\phi(z)=az+b$ with $|a|<1$, then
\beq\label{eq:conorm}
\|C_\phi\| = \exp \left( \frac14 \frac{|b|^2}{1-|a|^2} \right).
\eeq

Using these results we may prove the following.

\begin{thm}\label{th:asymp}
	The asymptotics of iteration of bounded composition operators on the Fock space $\GF^\nu$ are as follows.
\begin{enumerate}
\item If $\phi(z)=az$ with $|a|=1$ and $a \ne 1$ then the sequence $(C^n_\phi)$ consists of unitary operators; it
does not converge, even weakly.
\item If $\phi(z)=az+b$ with $|a|<1$ then the sequence $(C^n_\phi)$ converges in norm
to the operator mapping a function $f$ to the constant function with value $f(b/(1-a))$.
\end{enumerate}
\end{thm}

\beginpf
In the first case the operators $C^n_\phi$ are clearly unitary;
by considering the function $f(z)=z$ and $C^n_\phi f(z)=f(a^n z)$ we see that weak convergence does not
hold.

In the second case, 
the operator $C_\phi$ is power-bounded by \eqref{eq:conorm} since 
if $a_n=a^n$ and $b_n=\frac{1-a^n}{1-a}b$, then
$|b_n|^2/(1-|a_n|^2) \to |b|^2/|1-a|^2$.
Now we can use Proposition \ref{prop:2.1}. Since $C_\phi$ is compact, we have $r_e(C_\phi)=0$.
For the point spectrum we suppose there is an eigenvalue $\lambda$ 
of modulus 1 and eigenfunction $f$ satisfying $f(a z+b)=\lambda f(z)$.
Iterating, we have 
\beq\label{eq:itereig}
f(a^nz+ b(1-a^n)/(1-a)) = \lambda^n f(z).
\eeq
That is, $\lambda^n f(z)  \to f(b/(1-a))$; but $z$ was arbitrary, so $f$ is identically $0$ if $|\lambda|=1$ but $\lambda \ne 1$.

Finally, if $1 \in \sigma(C_\phi)$ then it is an eigenvalue, since $C_\phi$ is compact. Then \eqref{eq:itereig} still holds for an eigenfunction $f$, and
$f(z)=f(b/(1-a))$ so $f$ is constant.

Since $\ker (\Id-C_\phi)$ is   $1$-dimensional, we see that $1$ is a pole of the
resolvent of order at most 1.
Indeed, assume that the pole order of $1$ is larger than $1$. Looking at the Jordan normal form of $T_0:={C_\varphi}_{|X_0}$, where $X_0=P\FF^\nu$ and $P$ the residue, we see that there exists $f\in\FF^\nu$ such that $(\Id -C_\varphi)f=1_{\CC}$. Evaluating at $z=\frac{b}{1-a}$ we obtain a contradiction. Thus the pole order is $1$. It follows that $P$ is the projection onto $\ker (C_\varphi -\Id)=\CC 1_{\CC}$ along 
$\{f\in \FF^\nu:f(b/(1-a))=0\}$. Thus $Pf=f(b/(1-a))1_{\CC}$.     
This completes the proof that $C^n_\phi \to P$ in norm thanks to Proposition~\ref{prop:2.1}. 
Thus $(C^n_\phi)$ converges in norm, and the limit is the same as the pointwise limit, namely the rank-1 operator $T$
with $Tf  (z) = f(b/(1-a))$ for all $f$.
\endpf

\subsection{Semigroups}

The theory of $C_0$ semigroups (strongly continuous semigroups) is well described in many books, for
example \cite{applebaum,EN00,pazy}.

Recall that a $C_0$ semigroup $(T_t)_{t \ge 0}$ of bounded operators on a Banach space $X$ satisfies the conditions
\begin{enumerate}
\item $T_0=\Id$;
\item $T_{s+t}=T_s T_t$ for all $s, t \ge 0$;
\item the mapping $t \mapsto T_t x $ is continuous for all $x \in X$.
\end{enumerate}

The infinitesimal generator $A$ of such a semigroup $(T_t)_{t \ge 0}$  is defined by
\beq\label{eq:infgen}
Ax = \lim_{t \to 0}\frac{T_t x - x}{t} \qquad (x \in \dom(A)).
\eeq

In the situation when $T_t = W_{w_t,\phi_t}$ for all $t \ge 0$ we therefore have
\[
W_{w_{s+t},\phi_{s+t}}= W_{w_s,\phi_s} W_{w_t,\phi_t}= W_{w_s (w_t\circ \phi_s), \phi_t \circ \phi_s},
\]
so that 
\beq\label{eq:refws}
w_{s+t}=w_s(w_t \circ \phi_s)\mbox{  and }\phi_{s+t}=\phi_t \circ \phi_s.
\eeq

We write $\CC_-:=\{z \in \CC: \re z < 0\}$.

From now on all semigroups will be assumed to be strongly continuous.

\begin{thm}  \label{thm:compsg}
A semigroup of bounded composition operators $(T_t)_{t \ge 0}=  (C_{\phi_t})_{t \ge 0}$ on  ${\mathcal F}^\nu$
satisfies one of the following  for $t>0$:
\begin{enumerate}
\item $\phi_t(z)= \exp(\lambda t)z+C (\exp (\lambda t)-1) $ for some $\lambda \in \CC_- $ and $C \in \CC$;
\item $\phi_t(z)=\exp(\lambda t)z$ for some $\lambda \in i\RR$.
\end{enumerate}
The infinitesimal generator is given by $Af(z)=   \lambda(z+C)f'(z)$, where $C=0$ in case (2). Moreover  
$(T_t)_{t\geq 0}$ converges in norm as $t\to\infty$ if and only if $\lambda\in \CC_-$, and it converges 
to the operator mapping a function $f$ to the constant function with value $f(-C)$.
\end{thm}

\beginpf
From Proposition \ref{prop:11} we may write
 $\phi_t(z)=a_t z  + b_t $; we have $a_{t+s}=a_t a_s$, so $a_t=\exp(\lambda t)$
for some $\lambda \in \CC_-:=\{z \in \CC: \re z < 0\}$ (if $a_t$ is identically zero then the semigroup is not strongly continuous at $0$
since it would require $f(b_t) \to f(z)$ as $t \to 0$ for all functions $f$ and for all $z \in \CC$).
Also, $b_{s+t}=a_t b_s + b_t = a_s b_t + b_s$.
Thus  we have
\[
b_s (\exp (\lambda t) - 1)= b_t (\exp (\lambda s)-1),
\]
 i.e., $b_t= C (\exp (\lambda t)-1)$
for some constant $C \in \CC$.
If $\re \lambda = 0$, then by Proposition \ref{prop:wcobdd} we have boundedness only for $C=0$.

Then we have
\[
Af(z)= \lim_{t \to 0} \frac{f(e^{\lambda t}z + C(e^{\lambda t}-1))-f(z)}{t}= \lambda(z+C)f'(z).
\]
The last assertion is an obvious corollary of  Theorem~\ref{th:asymp}.
\endpf

A semigroup $(T_t)_{t \ge 0}$ is said to be quasicontractive if
for some $\omega \in \RR$
$\|T_t\| \le e^{\omega t}$ for all $t \ge 0$.
 
\begin{rem}
If $\phi_t(z)= \exp(\lambda t)+C (\exp (\lambda t)-1) $ with $\re\lambda<0$ then
by \eqref{eq:conorm}, on the Fock space $\GF^\nu$, we have
\beq 
\|T_t\| = \exp \left( \frac14 \frac{|C^2  (\exp (\lambda t)-1)|^2}{1-|e^{2\lambda t}|} \right).
\eeq

A routine calculation shows that
\[
 \frac{|\exp (\lambda t)-1|^2}{1-|e^{2\lambda t}|} = O(t) 
 \]
 as $t \to 0$. Hence $(T_t)$ is a quasicontractive semigroup.
 \end{rem}

In \cite[Thm.~3.9]{acp15} It was shown that the following conditions are equivalent for a semigroup with 
generator $A: f \mapsto G(z)f'(z)$ on the Hardy space $H^2$, where $G \in H^2$:
\begin{enumerate}
\item $A$ generates a $C_0$ semigroup of composition operators on $H^2(\DD)$;
\item $A$ generates a quasicontractive semigroup on $H^2(\DD)$;
\item $\esssup_{z \in \TT} \re (\overline z G(z)) \le 0$.
\end{enumerate}

In the light of Theorem \ref{thm:compsg} and the following proposition we see that similar condition on $G$  applies in the case of the Fock space, namely
\beq\label{eq:limsup}
\limsup_{|z| \to \infty} \re (\overline z G(z)) \le 0.
\eeq

\begin{prop}
For an entire function $G$, condition \eqref{eq:limsup} holds if and only if   $G(z)=a z+b$ with either (i)~$a \in \CC_-$ and $b \in \CC$,
or else (ii)~$a \in i\RR$ and $b=0$. 
\end{prop}

\beginpf
Since $\overline z (az+b)= a|z|^2 + b \overline z$, it is clear that functions $G$ of the form above do satisfy \eqref{eq:limsup}.

Conversely, for an entire function $G$  let $F(z)=(G(z)-G(0))/z$. Thus $\overline z G(z) = |z|^2 F(z) + \overline z G(0) $.  For $G$ to satisfy 
\eqref{eq:limsup} we must have $\re (|z|^2 F(z)+ \overline z G(0)) < 1 $, and hence
$\re   (F(z)+   G(0)/z) < 1 $, for sufficiently large $|z|$. This is clearly impossible if $F$ is a nonconstant polynomial, and also impossible if $F$
has an essential singularity at $\infty$ (for example, by the Weierstrass--Casorati theorem). So $F$ must be constant, and the conditions on $a$ and $b$ 
are now immediate.
\endpf

Next we consider the numerical range of a semigroup generator $A: f \mapsto G(z)f'(z)$ on the Fock space $\GF^2$,
namely
\[
W(A) = \{ \langle Af,f \rangle: f \in \GD(A) , \|f\|=1\},
\]
where $\GD(A)$ denotes the domain of $A$. Recall that a semigroup is quasicontractive if and only if $W(A)$ is contained in some left half-plane.

\begin{prop}
If an entire function $G$ is such that the operator $A: f \mapsto Gf'$ is the generator of a quasicontractive $C_0$ semigroup on $\GF^2$, then
$G$ is a polynomial of degree at most $1$, say $G(z)=az+b$, with $\re a \le 0$.
\end{prop}
\beginpf
 We write $G(z)=\sum_{k=0}^\infty g_k z^k$ and examine the coefficients $g_k$.
First, we calculate $\langle Af,f \rangle$ for $f(z)=e_n(z)$ as defined in \eqref{eq:onb}.
\[
\langle Af,f \rangle = \langle Gf',f \rangle = g_1 n \langle z^n,z^n \rangle / n! = ng_1.
\]
Now if $A$   generates a quasicontractive $C_0$ semigroup, there must be an upper bound on $\re W(A)$, and hence $\re g_1 \le 0$.

Next, consider the function $f(z) = a z^n + b z^{n+m}$ where $m \ge 1$ is fixed. We see that 
$\|f\|^2 = |a|^2 n! + |b|^2 (n+m)!$.

Also
\begin{eqnarray*}
\langle Af,f \rangle &=& \langle Gf',f \rangle\\
& =& g_1(n|a|^2 n! + (n+m)|b|^2 (n+m)! ) + g_{m+1} (n+m)a \overline b (n+m)!.
\end{eqnarray*}
Now take $a=n^{m/2}\omega$ where $|\omega|=1$, and $b=1$, so that $\|f\|^2 = O(n! n^m)$ and
\[
\langle Af,f \rangle =  g_1(n^{m+1} n! +  (n+m) (n+m)! ) + g_{m+1}(n+m)(n^{m/2} \omega)   (n+m)!.
\]
That is, if $g_{m+1} \ne 0$, then for a suitable choice of $\omega$,  $\re \langle Af,f \rangle / \|f\|^2$ grows as $|g_{m+1}| n^{m/2+1}$ as $n \to \infty$;
hence $g_{m+1}=0$ for all $m \ge 1$ if $A$ is the generator of a quasicontractive $C_0$ semigroup.
\endpf

We can now go further and describe all semigroups with generators $A:f \mapsto Gf'$.
In \cite{AC19} there is a result applying to function spaces, i.e., Banach spaces $X$ of holomorphic functions on a domain $\Omega$ for which
point evaluations are continuous functionals. The result applies to such spaces satisfying the following supplementary
condition:

(E) \quad 
If $(z_n)$ is a sequence in $\Omega$ such that $z_n \to z \in \overline\Omega \cup \{\infty\}$ and $\lim f(z_n)$ exists in $\CC$  for all $f \in X$
then $z \in \Omega$.

\begin{thm}\label{thm:E}
Suppose that $(T_t)_{t \ge 0}$ is a $C_0$ semigroup on a function space $X$ with property (E) such that the generator $A$
is the operator $f \mapsto Gf'$ for all $f \in \dom(A)$ for some $G \in \Hol(\Omega)$. Then $(T_t)$ is a semigroup of composition operators.
\end{thm}

Similar results are given in \cite{GY19}, for function spaces on the disc $\DD$.\\

Theorem \ref{thm:E} clearly applies to the space $X=\GF^2$ with $\Omega=\CC$ and so we have the following corollary.

\begin{cor}
Every $C_0$  semigroup   on $\GF^2$ with generator $A: f \mapsto Gf'$ for some $G \in \Hol(\CC)$ 
consists entirely of composition operators.
\end{cor}

We may also consider the case of analytic semigroups, analogously to \cite{acp16}. For $0< \theta < \pi/2$ we define the sector $\Sigma_\theta$ by
\[
\Sigma_\theta = \{re^{i\alpha}: r>0, |\alpha|<\theta \},
\]
and recall that an analytic (or holomorphic) semigroup is one for which there is an analytic mapping $T: \Sigma_\theta \to \LL(X)$
satisfying the semigroup property,  and  which, on defining $T(0)=\Id$, is strongly continuous at $0$, and satisfies $\sup_{\xi \in \Sigma_\theta \cap \DD} \|T_\xi\| < \infty$.

We may test   the analytic property by using Corollary 2.3 of \cite{acp16}, which asserts that if $A$ is a Hilbert-space operator and
$A$ and $e^{\pm i\theta}A$ generate quasicontractive semigroups for some $0 < \theta < \pi/2$, then $A$ generates an analytic
semigroup on the sector $\Sigma_\theta$.

As a variation on Corollary 2.5 of \cite{acp16}, we also have the following result.
\begin{lem}
Let $(T_t)_{t \ge 0}$ be a semigroup of composition operators on $\GF^2$ and suppose that it has an analytic extension
to $\Sigma_\theta$ for some $0<\theta < \pi/2$. Then for every $\xi \in \Sigma_\theta$ the operator $T_\xi$ is 
a composition operator.
\end{lem}
\beginpf
We note first that, writing $e_n(z)=z^n$ for $z \in \CC$ we have $T_t(e_n) = (T_t(e_1))^n$ for all $n \in \NN$ and $t \ge 0$.
Now defining $F(\xi)=T_\xi e_n - (T_\xi e_1)^n$ we have that $\xi \mapsto F(\xi)(z_0)$ is an
analytic function on $\Sigma_\theta$ for each $z_0 \in \CC$; it vanishes on $\RR_+$, and hence is identically zero.  Thus by
linearity and continuity $T_\xi=C_{\phi_\xi}$ for each $\xi$, where $\phi_\xi = T_\xi e_1$.
\endpf

We thus have the following result. 
\begin{thm}\label{thm:analytic}
For an entire function $G$   the operator $A:f \mapsto Gf'$ is the generator
of an analytic semigroup of operators on $\GF^2$ if and only if there is $\theta \in (0,\pi/2)$ such that
\beq 
\limsup_{|z| \to \infty} \re   (e^{\pm i\theta} \overline z G(z)) \le 0.
\eeq
In this case, the semigroup consists of composition operators and is quasicontractive.
\end{thm}

\section{Weighted composition operators}\label{sec:3}

As is the case for weighted composition operators on other function spaces, the
results that one can prove are necessarily less complete than the results for
(unweighted) composition operators. Nevertheless we may make some progress here.

\subsection{Iteration}


For the weighted composition operator $\W$ with $\phi(z)=az+b$ and $a \ne 1$ we have
\[
\W^n f(z)= w(z)w(\phi(z))w(\phi^2(z)) \ldots w(\phi^{n-1}(z))
f\left(a^n z + \frac{1-a^n}{1-a}b\right)
\]
for $n=1,2,3,\ldots$, where we denote the $k$th iterate of $\phi$ by $\phi^k$.

For $a=1$ the corresponding formula is
\[
\W^n f(z)= w(z)w(\phi(z))w(\phi^2(z)) \ldots w(\phi^{n-1}(z))
f\left(  z + nb\right).
\]

A necessary and sufficient condition for boundedness of $\W$ on $\GF^\nu$ is given above in
Proposition  \ref{prop:wcobdd}.



We do not have an exact expression for the norm of a weighted composition operator, but
from \cite{HaiKhoi} we know that
for $\phi(z)=az+b$ with $|a| \le 1$
we have
\beq\label{eq:hokey}
 \sqrt{M(w,\phi)} \le \|W_{w,\phi}\| \le  \sqrt{M(w,\phi)} /|a|.
 \eeq
Writing $M_z(w,\phi)= |w(z)|^2 e^{(|\phi(z)|^2-|z|^2)}$ we have that $\W$ is compact
if and only if $|a|<1$ and $M_z(w,\phi) \to 0$ as $|z| \to \infty$.\\

For power-boundedness, some results are given in \cite{Men22}, with applications to mean ergodicity.
We proceed independently, and our main concern is with discrete and continuous semigroups.

First, it is possible to give an exact answer in the case $|a|=1$.

\begin{thm}\label{thm:3.1}
Suppose that $\phi(z)=az+b$ with $|a|=1$. Then\\
$\W$ is power-bounded on $\GF^\nu$ if and only if $|w(0)| \le e^{-|b|^2/2}$.
\end{thm}

\beginpf
(i) We begin with the case $a=1$, and with $w(z)=w(0)e^{-\overline b a z}$ we have 
\[
\W^n f(z)=w(z)w(z+b)\ldots w(z+(n-1)b)f(z+nb).
\]
That is, it is a weighted composition operator of the form $W_{w_n,\phi_n}$, and we have
\begin{eqnarray*}
M_z(w_n,\phi_n) &=& |w_n(z)|^2 e^{|\phi_n(z)|^2-|z|^2} \\
&= &
|w(0)|^{2n}|e^{-2\overline b z}| \, |   e^{-2\overline b(z+b)}| \ldots |e^{-2\overline b(z+(n-1)b)}|
e^{|z+nb|^2-|z|^2} \\
&=& |w(0)|^{2n} |e^{-2n\overline b z}| e^{-n(n-1)|b|^2} e^{2\Re nz\overline b}e^{n^2|b|^2} \\
&=& |w(0)|^{2n} e^{n|b|^2} 
\end{eqnarray*}
and this is uniformly bounded in $n$ (and $z$) if and only if $|w(0)| \le e^{-|b|^2/2}$.\\
(ii) For $a \in \TT \setminus\{1\}$ we have
\begin{eqnarray*}
\W^n f(z)&=& w(z)w(az+b)w\left(a^2z+\frac{1-a^2}{1-a}b\right) \ldots w\left(a^{n-1}z+\frac{1-a^{n-1}}{1-a}b\right) \\
&& \times 
f\left(a^n z + \frac{1-a^n}{1-a}b\right),
\end{eqnarray*}
which again we write as $W_{w_n,\phi_n}$. We now have
\[
w_n(z)= w(0)^n  e^{-\overline b a z}e^{-\overline b a (az+b)}\ldots e^{-\overline b a (a^{n-1}z+\frac{1-a^{n-1}}{1-a}b)}
\] 
and
\begin{eqnarray} \label{eq:mzn}
M_z(w_n,\phi_n) &=& |w_n(z)|^2 e^{|\phi_n(z)|^2-|z|^2} \nonumber\\
&= & |w(0)|^{2n} e^{-2|b|^2 \re (a(1+\frac{1-a^2}{1-a}+ \ldots + \frac{1-a^{n-1}}{1-a}))}
e^{-2 \re( \overline b a z\frac{1-a^n}{1-a})}\nonumber \\ && \times 
e^{|a^n z + \frac{1-a^n}{1-a}b|^2-|z|^2}
\end{eqnarray}
and
\[
e^{|a^n z + \frac{1-a^n}{1-a}b|^2-|z|^2} = e^{|\frac{1-a^n}{1-a}b|^2 + 2 \re (a^n z\overline b \frac{1-\overline a^n}{1-\overline a})} = 
e^{|\frac{1-a^n}{1-a}b|^2 + 2 \re (a z\overline b \frac{1-  a^n}{1-  a})} 
\]
since $|a|=1$.

This gives
\[
M_z(w_n,\phi_n) =  |w(0)|^{2n} e^{|\frac{1-a^n}{1-a}b|^2}
e^{-2|b|^2 \re [ (n-1)a/(1-a)- a^2(1-a^{n-1})/(1-a)^2 ]}
\]
Note that for $|a|=1$ and $a \ne 1$ we have 
\[
\re(a/(1-a))=\frac12\left(\frac{a}{1-a}+\frac{\overline a}{1-\overline a}\right)=
\frac12\left(\frac{a}{1-a}+\frac{1}{a-1}\right)=-\frac12,
\]
so that 
\[
\sup_{z \in \CC} M_z(w_n,\phi_n) \le K |w(0)|^{2n} e^{|b|^2(n-1)},
\]
where $K$ is independent of $n$;  hence uniform boundedness of
the powers $\W^n$ on $\GF^\nu$ is again equivalent to the
condition  $|w(0)| \le e^{-|b|^2/2}$.
\endpf

The analysis is apparently intractable in the case $|a|<1$, however we may consider the case when the weight $w$ does
not vanish, which will be of importance when we come to consider semigroups. In this case, writing $\beta=1-|a|^2$ we have
the following result from \cite{CG21}.
\begin{thm}\label{thm:3.2}
For $\phi(z)=az+b$ with $|a|<1$ and 
$w$ nonvanishing the operator $\W$ is bounded on $\GF^\nu$ if and only if
\[
w(z)= e^{p+qz+rz^2}
\]
and either:\\
(i) $|r|< \beta/2$, in which case $\W$ is compact on $\GF^\nu$, or\\
(ii) $|r|=\beta/2$ and, with $t=q+\overline b a$, one has either $t=0$ or else
$t \ne 0$ and $r=-\frac{\beta}2\frac{t^2}{|t|^2}$. In case (ii) $\W$ is not compact on $\GF^\nu$.
\end{thm}


The analysis splits into two cases:

{\bf Case 1.} We take $\phi(z)=az+b$ with $0<|a|<1$, and $w(z)=e^{p+q z}$. That is, $r=0$.

 Here we once again have \eqref{eq:mzn} with $w(0)=e^{p}$ but to check boundedness all we have
 is \eqref{eq:hokey}.
 
{\bf Case 2.} We  take $\phi(z)=az+b$ with $0<|a|<1$, and $w(z)=e^{p+qz+ r z^2}$
with $0 < |r| \le \beta/2$.

In Case 1, we now have
\begin{eqnarray}  \label{eq:case1}
M_z(w_n,\phi_n) &=& |w_n(z)|^2 e^{|\phi_n(z)|^2-|z|^2} \nonumber\\
&=& |w(0)|^{2n} |e^{2q(z+az+\ldots+a^{n-1}z)}| | e^{2q b (1+\frac{1-a^2}{1-a}+ \ldots \frac{1-a^{n-1}}{1-a})}|
e^{|\phi_n(z)|^2-|z|^2} \nonumber\\
&= & e^{2n \re p}
e^{2\re [q z   (1-a^{n})/(1-a)  ]}
e^{2 \re q b[n  /(1-a)-a(1-a^{n-1})/(1-a)^2]}\nonumber \\ && \times 
e^{|a^n z + \frac{1-a^n}{1-a}b|^2-|z|^2} 
\end{eqnarray}

Taking logarithms,  we 
have the following:
\begin{thm}\label{thm:3.3}
In Case 1:\\
if $\re p +  \re \left( \frac{ q b}{1-a}\right) > 0$ then $\W$ is not power-bounded on $\GF^\nu$;\\
if $\re p +   \re \left( \frac{ q b}{1-a}\right) - \ln |a| < 0 $ then $\W$ is power-bounded on $\GF^\nu$.
\end{thm}
\beginpf
In \eqref{eq:case1} we have that $a^n \to 0$, and the asymptotic behaviour of
$\sup_{z} M_z(w_n,\phi_n)$ depends only on the first and third factors. Now use 
\eqref{eq:hokey}.
\endpf

Note that this matches the condition of Theorem \ref{thm:3.1} if $|a|=1$ and $q=-\overline ba$.

In Case 2 we have the extra factor $v(z)=e^{rz^2}$ in the weight, and

Then
\begin{eqnarray*}
v(z)v(\phi(z))\ldots v(\phi_{n-1}(z))
&=& e^{r z^2}
e^{r(az+b)^2}\ldots 
e^{r(a^{n-1}z+ b \frac{1-a^{n-1}}{1-a})^2} \\
&=& e^{A+Bz+Cz^2} ,
\end{eqnarray*}
where 
\begin{eqnarray*}
A &=& r b^2 \left(\left(\frac{1-a}{1-a} \right)^2+ \left(\frac{1-a^2}{1-a} \right)^2 + \ldots + \left(\frac{1-a^{n-1}}{1-a}\right )^2\right)\\
&=& \frac{r b^2}{(1-a)^2}\left ( n-1 - 2a \frac{1-a^{n-1}}{1-a} + a^2 \frac{1-a^{2n-2}}{1-a^2}
\right),\\
B &=& 2rb \left (a\frac{1-a}{1-a}+\ldots+a^{n-1}\frac{1-a^{n-1}}{1-a} \right ) \\
&=&
\frac{2rb}{1-a} \left( a \frac{1-a^{n-1}}{1-a}- a^2 \frac{1-a^{2n-2}}{1-a^2} \right )\\
\noalign{\hbox{and}}\\
C &=& r+r a^2 + \ldots r a^{2n-2} = r \frac{1-a^{2n}}{1-a^2}.
\end{eqnarray*}
Recall that $|r| \le (1-|a|^2)/2$. 

We therefore have the following result:
\begin{thm}\label{thm:3.4}
In Case 2:\\
if $\re p +  \re \left( \frac{ q b}{1-a}\right) + \re \left( \frac{r b^2}{(1-a)^2} \right)> 0$ then $\W$ is not power-bounded on $\GF^\nu$;\\
if $\re p +   \re \left( \frac{ q b}{1-a}\right) + \re \left( \frac{r b^2}{(1-a)^2} \right) - \ln |a| < 0 $ 
then $\W$ is power-bounded on $\GF^\nu$.
\end{thm}
\beginpf
We again have that $a^n \to 0$, and the asymptotic behaviour of
$\sup_{z} M_z(w_n,\phi_n)$ is as in Case 1, with the additional terms from A and C. Once more we use 
\eqref{eq:hokey}.
\endpf

Coming to the question of convergence of iterates, we note that if $|r| < (1-|a|^2)/2$, then $\W$ is compact,
and so its essential spectrum is $\{0\}$. We now consider the point spectrum.

\begin{thm}\label{th:eigen}
For a compact weighted composition operator $\W$ on $\GF^\nu$ with $w$ not the zero function, $\phi(z)=az+b$ and $0< |a| < 1$, the point spectrum
consists of simple eigenvalues.
\end{thm}
\beginpf
The mapping $\phi$ has a Denjoy--Wolff point $\alpha \in \CC$, i.e., a solution to $\phi(\alpha)=\alpha$, namely 
$\alpha=b/(1-a)$, to which all its iterates converge pointwise.

Factorizing out any possible zeros at $z=0$, we may write 
\[
w(\alpha+z)= z^N \exp (\sum_{n=0}^\infty d_n z^n),
\]
and let us suppose that $\W f=\lambda f$ with 
\[
f(\alpha+z)=z^M \exp (\sum_{n=0}^\infty c_n z^n).
\]
Note that these expansions are valid in a disc centered at $z=0$.

Now $\phi(\alpha+z)=a(\alpha+z)+b=\alpha+az$, so
we have  
\[
z^N \exp (\sum_{n=0}^\infty d_n z^n) a^M z^M \exp (\sum_{n=0}^\infty  c_n a^n z^n) = \lambda z^M \exp (\sum_{n=0}^\infty c_n z^n).
\]
We thus have $N=0$, and $\lambda=a^M e^{d_0}=a^M w(\alpha)$; moreover, for some $k \in \ZZ$,
\[
 \sum_{n=1}^\infty  (d_n+c_n a^n) z^n = \sum_{n=1}^\infty c_n z^n + 2\pi i k\]
at least on some disc   $\{|z|< r \}$. So $k=0$ and $c_n=d_n/(1-a^n)$.

Thus the eigenvector (if it is in the space) is unique up to a constant factor.

\endpf

Putting all these ideas together along with Proposition \ref{prop:2.1} we arrive at the following result.

\begin{cor}\label{cor:cv}
For a  nontrivial compact and power-bounded  weighted composition $\W$ on $\GF^\nu$ with $\phi(z)=az+b$ and $0<|a|<1$,
let $\alpha=b/(1-a)$. Then,
provided that $a^M w(\alpha) \not\in \TT\setminus \{1\}$ for all $M=0,1,2,\ldots$, 
 the sequence $(\W^n)$
converges in norm to the operator mapping $f$ either to the constant function $f(\alpha)$ if $w(\alpha)=1$, otherwise to the zero
function.
\end{cor}

\begin{rem}
Alternative approaches to Theorems \ref{thm:3.1}, \ref{thm:3.3}, \ref{thm:3.4} and \ref{th:eigen} have been suggested by a referee.
Namely, in Theorem \ref{thm:3.1} one can show that the operator $\W$ can be written as $w(0)e^{|b|^2/2}$
times a surjective isometry on $\GF^\nu$. In Theorems \ref{thm:3.3}, \ref{thm:3.4} and \ref{th:eigen} the operator
$\W$ is similar to an operator of the form $W_{\tilde w, \tilde\phi}$ where $\tilde \phi(z)=az$: this approach may simplify the calculations.
\end{rem}

\subsection{Semigroups and groups}

For semigroups of weighted composition operators, say, $T_t f(z)= w_t(z) f(\phi_t(z))$, we  have the relation
\eqref{eq:refws}
so if $\phi_t(z)=a_tz+b_t$ we obtain
\beq\label{eq:wrule}
w_{t+s}(z)= w_s(z) w_t(a_s z+b_s) = w_t(z) w_s(a_t z + b_t)
\eeq
while
$a_{t+s}=a_t a_s$ and 
$b_{s+t}=a_t b_s + b_t = a_s b_t + b_s$.

\begin{thm}\label{thm:c0sg}
A $C_0$ semigroup $T_t f(z)= w_t(z) f(\phi_t(z))$ ($t \ge 0$) of weighted composition operators on the Fock space  $\GF^\nu$ has one of the
following  two expressions;
\begin{enumerate}
\item $\phi_t(z)= \exp(\lambda t)z+C (\exp (\lambda t)-1) $ for some $\lambda \in \CC_- $ and $C \in \CC$;
in which case $w_t=e^{p_t+q_t z+r_t z^2}$ with $p_t,q_t,r_t$ given respectively by \eqref{eq:pt}, \eqref{eq:qt} and \eqref{eq:rt} below.
\item $\phi_t(z)=\exp(\lambda t)z+C (\exp (\lambda t)-1) $ for some $\lambda \in i\RR$ and $C \in \CC$, in which case 
$w_t(z)=w_t(0)\exp(\overline C (\exp (\lambda t)-1)  z)$
and moreover 
$w_t(0)=e^{\mu t}e^{|C|^2 (e^{\lambda t}-1)}$ for some $\mu \in \CC$.
\end{enumerate}
\end{thm}

\beginpf
The formula for $\phi_t$ follows as in Theorem \ref{thm:compsg}. We now have two cases to analyse:

{\bf (1).} The proof of \cite[Lem. 2.1(b)]{koenig}, can be adapted and simplified to show that $w_t(z) \ne 0$ for all $t \ge 0$ and $z \in \CC$.

For if there are a $z_0 \in \CC$ and $t_0\ge 0$ with $w_{t_0}(z_0)=0$, then let $I=\{t \ge 0: w_t(z_0)=0 \}$ and $\tau=\inf I$.
Since the semigroup is strongly continuous we have $\tau \in I$ 
by considering the effect of the semigroup on a constant function; also $\tau>0$ since $w_0(z)=1$ for all $z$.
By \eqref{eq:wrule} we have $I=[\tau,\infty)$ and, taking $s=\tau/2$, we have 
\[
w_{t}(z_0)w_s(a_tz_0+b_t)=w_{t+s}(z_0)=0
\] 
for all $t \in (\tau/2,\tau)$
and hence
$w_s(a_tz_0+b_t)=0$
for all such $t$. Hence $w_s$ is identically zero, as it is zero on a nontrivial continuous path. \\

Now, by Theorem \ref{thm:3.2}, $w_t(z)$ takes the form 
\[
w_t(z)=\exp(p_t+q_tz+r_t z^2).
\]

From \eqref{eq:wrule} we have
\[
p_{t+s}+q_{t+s}z + r_{t+s}z^2 = p_s+q_s z + r_s z^2+p_t+q_t(a_s z+b_s)+r_t(a_sz + b_s)^2,
\]
For extra legibility, write $x=\exp(\lambda s)$ and $y=\exp(\lambda t)$; then 
$a_s=x$ and $b_s=C(x-1)$; so
the
equations we have are
\[
p_{t+s}+q_{t+s}z + r_{t+s}z^2 = p_s+q_s z + r_s z^2+p_t+q_t(a_s z+b_s)+r_t(a_sz + b_s)^2,
\]
giving
\begin{eqnarray*}
p_{t+s} &=& p_s + p_t + q_tb_s+r_tb_s^2 , \\
q_{t+s} &=& q_s+q_t a_s + 2r_t a_sb_s, \\
r_{t+s} &=& r_s + a_s^2 r_t.
\end{eqnarray*}
Thus $r_{t+s}=r_s+r_tx^2$, for which the general solution is 
$
r_s=\alpha(1-x^2)
$
or
\beq\label{eq:rt}
r_t=\alpha(1-e^{2\lambda t})
\eeq
 for some $\alpha \in \CC$.
(It is easy to check that $r_t/(1-e^{2\lambda t})$ is necessarily constant.)
The constant $\alpha$ must also be such that $|r_t| \le (1-|a_t|^2)/2$ by virtue of Theorem \ref{thm:3.2}.\\

Next $q_{t+s}=q_s + q_t x + 2\alpha C(1-y^2)x(x-1)$, which has solution
$
q_s = 2\alpha C x(1-x) + \gamma (1-x)
$
or
\beq\label{eq:qt}
q_t=2\alpha C e^{\lambda t}(1-e^{\lambda t}) + \gamma (1-e^{\lambda t})
\eeq
 for some $\gamma \in \CC$.\\

Finally,
$p_{t+s}=p_s+p_t +2\alpha C^2  y(1-y)(x-1) + C\gamma(1-y)(x-1)+\alpha C^2 (1-y^2) (x-1)^2$,
giving
$
p_s=-\alpha C^2 (x-1)^2 - \gamma C(x-1) + \delta s
$
or
\beq\label{eq:pt}
p_t=-\alpha C^2 (e^{\lambda t}-1)^2 - \gamma C(e^{\lambda t}-1) + \delta t
\eeq
 for some $\delta \in \CC$.\\

{\bf (2).}
We require $w_t(z)=w_t(0)e^{-\overline {b_t}a_t z}$ to satisfy \eqref{eq:wrule}, that is,
\[
w_{t+s}(0) e^{-\overline{b_{s+t}}a_{s+t}z} = w_s(0) e^{-\overline{b_s}a_s z}w_t(0) e^{-\overline{b_t}a_t (a_s z + b_s)}.
\]
Observe that
\[
\overline{b_{s+t}}a_{s+t} = \overline{b_s}a_s + \overline{b_t}a_t a_s ,
\]
that is,
\[
\overline C (e^{-\lambda (s+t)}-1)e^{\lambda (s+t)} = \overline C (e^{-\lambda s}-1)e^{\lambda s} + \overline C(e^{-\lambda t}-1)e^{\lambda t}e^{\lambda s},
\]
since
\[
1-e^{\lambda (s+t)} = 1-e^{\lambda s} + (1-e^{\lambda t})e^{\lambda s}
\]
holds for all $s$ and $t$.

We also have 
\[
w_{t+s}(0) = w_s(0)  w_t(0)e^{-\overline{b_t}a_t b_s} = w_t(0)  w_s(0)e^{-\overline{b_s}a_s b_t}
\]
since for this to happen we require
\[
|C|^2 (1-e^{\lambda t} )(e^{\lambda s}-1) = |C|^2 (1-e^{\lambda s} )(e^{\lambda t}-1) 
\]
which holds for all $s$ and $t$. Finally,
\[
w_{t+s}(0) = w_s(0)  w_t(0) e^{|C|^2 (e^{\lambda t}-1 )(e^{\lambda s}-1)}
\]
or,
with $v_t=w_t(0) e^{-|C|^2 (e^{\lambda t}-1)}$,
\[
v_{t+s}=v_s  v_t e^{|C|^2((e^{\lambda t}-1 )(e^{\lambda s}-1)-(e^{\lambda (t+s)}-1)+(e^{\lambda t}-1) + (e^{\lambda s}-1)}=v_s v_t
\]
so that finally $w_t(0)=e^{\mu t}e^{|C|^2 (e^{\lambda t}-1)}$ for some $\mu \in \CC$.
\endpf

\begin{rem}We may also calculate   infinitesimal generators 
using \eqref{eq:infgen}, so that for $T_t(f)(z)=w_t(z)f(\phi_t(z))$
with \[w_t(z)=\exp(p_t+q_t z+r_t z^2) \mbox{ and }\phi_t(z)=a_tz+b_t\]
we have
\begin{eqnarray*}
Af(z) &=&  \lim_{t \to 0}w_t(z) \frac{f \circ \phi_t(z)-f(z)}{t} + \lim_{t \to 0}\frac{w_t f(z)-f(z)}{t}\\
&=& \left[\frac{\partial}{\partial t}(a_t z+ b_t ) f'(z)+ \frac{\partial}{\partial t}(  p_t +   q_t z + r_t z^2) f(z)\right]_{t=0}.
\end{eqnarray*}
\end{rem}


Finally, we may characterise one-parameter groups of weighted composition operators, using  Proposition \ref{prop:invert} and Theorem \ref{thm:c0sg}.

\begin{cor}
A $C_0$ group $T_t f(z)= w_t(z) f(\phi_t(z))$ ($t \in \RR$) of weighted composition operators on the Fock space $\GF^\nu$ satisfies
\[\phi_t(z)=\exp(\lambda t)z+C (\exp (\lambda t)-1)\]
 for some $\lambda \in i\RR$ and $C \in \CC$, in which case 
\[w_t(z)=w_t(0)\exp(\overline C (\exp (\lambda t)-1)  z)\]
and moreover 
$w_t(0)=e^{\mu t}e^{|C|^2 (e^{\lambda t}-1)}$ for some $\mu \in \CC$.
\end{cor}
\beginpf
The first case of 
 Theorem \ref{thm:c0sg} does not give bounded inverses, but
 in the second case where $a_t=e^{\lambda t}$
and $b_t=C(e^{\lambda t}-1)$ with $\lambda \in i\RR$,
we do have $a_{-t}=a_t^{-1}$ and $b_{-t}=-a_t^{-1}b_t$; also 
$w_{-t}(0)=w_t(0)^{-1}e^{-|b_t|^2}$, so that $T_{-t}=T_t^{-1}$ by
Proposition  \ref{prop:invert}.
\endpf

\begin{rem}
It is also possible to consider the question of when a weighted composition semigroup has an
analytic extension, analogously to Theorem \ref{thm:analytic}. However, the question is
simply whether $e^{\pm i\theta}A$ also generates a semigroup for some $0 < \theta < \pi/2$, which reduces to
considering the semigroups in Theorem \ref{thm:c0sg} with a complex parameter. Accordingly,
we omit the details.
\end{rem}

\subsection*{Acknowledgment}

The authors are grateful to a referee for a very careful of the reading of the paper and
the correction of numerous errors.


\end{document}